\documentclass[12pt]{article}
\usepackage{amssymb,amsmath}
\usepackage{cases}
\usepackage{amsfonts}
\usepackage{cite,color,xcolor}
\usepackage[left=2.0cm,right=2.0cm,top=2.0cm,bottom=2.0cm]{geometry}
\usepackage[colorlinks,citecolor=blue,urlcolor=blue]{hyperref}

\newtheorem{theorem}{Theorem}[section]

\newtheorem{lemma}{Lemma}[section]
\newtheorem{proposition}{Proposition}[section]

\newtheorem{remark}{Remark}[section]

\newcommand{\bal}{\begin{align}}
\newcommand{\bbal}{\begin{align*}}
\newcommand{\beq}{\begin{equation}}
\newcommand{\eeq}{\end{equation}}
\newcommand{\bca}{\begin{cases}}
\newcommand{\eca}{\end{cases}}
\newcommand{\pa}{\partial}
\newcommand{\fr}{\frac}
\newcommand{\na}{\nabla}
\newcommand{\De}{\Delta}

\newcommand{\cd}{\cdot}
\newcommand{\ep}{\varepsilon}
\newcommand{\dd}{\mathrm{d}}

\newcommand{\R}{\mathbb{R}}

\newcommand{\les}{\lesssim}

\newcommand{\D}{\mathrm{div}}

\newcommand{\bi}{\Big}

\begin{document}
\title{Non-uniform dependence on initial data for the 2D viscous shallow water equations}

\author{Jinlu Li$^{1}$, Yanghai Yu$^{2,}$\footnote{E-mail: lijinlu@gnnu.edu.cn; yuyanghai214@sina.com(Corresponding author); mathzwp2010@163.com} and Weipeng Zhu$^{3}$\\
\small $^1$ School of Mathematics and Computer Sciences, Gannan Normal University, Ganzhou 341000, China\\
\small $^2$ School of Mathematics and Statistics, Anhui Normal University, Wuhu 241002, China\\
\small $^3$ School of Mathematics and Information Science, Guangzhou University, Guangzhou 510006, China}

\date{\today}

\maketitle\noindent{\hrulefill}

{\bf Abstract:} The failure of uniform dependence on the data is an interesting property of classical solution for a hyperbolic system. In this paper, we consider the solution map of the Cauchy problem to the 2D viscous shallow water equations which is a hyperbolic-parabolic system. We prove that the solution map of this problem is not uniformly continuous in Sobolev spaces $H^s\times H^{s}$ for $s>2$.


\vskip0mm\noindent{\hrulefill}

\section{Introduction}\label{sec1}
The 2D viscous shallow water equations are given by the following system, which have systematically introduced in \cite{Bresch1,Bresch2}
\begin{eqnarray}\label{o}
        \left\{\begin{array}{ll}
          \partial_t\rho+\D (\rho u)=0,\\
          \partial_t(\rho u)+\D(\rho u\otimes u)-\mu\D(\rho\nabla u)+\rho\na \rho=0,\\
           \rho|_{t=0}=\rho_0,\;u|_{t=0}=u_0.
          \end{array}\right.
        \end{eqnarray}
where $\rho(x,t)$ is the height of fluid surface, $u(x,t)$ is the horizontal velocity field and $\mu>0$ is the viscous coefficient. We suppose that
the initial data $\rho_0(x)$ is a small perturbation of some positive constant $\overline{\rho}_0$.

Classical solutions and well-posedness of the initial (boundary) value problem for the shallow water equations \eqref{o} have been studied extensively. By using Lagrangian coordinates and H\"{o}lder space estimates, Bui \cite{Bui} obtained the local existence and uniqueness of classical solutions to the Cauchy-Dirichlet problem for \eqref{o} with initial data in $C^{2+\alpha}$. With the help of the energy method of Matsumura and Nishida \cite{Matsumura}, Kloeden \cite{Kloeden} and Sundbye \cite{Sondbye1} independently showed the global existence and uniqueness of
classical solutions to the Cauchy-Dirichlet problem for \eqref{o}. Subsequently, Sundbye \cite{Sondbye2} also proved the
existence and uniqueness of classical solutions to the Cauchy problem for \eqref{o} using the method of \cite{Matsumura}. Due to the strong nonlinearity of system \eqref{o}, the problem of existence of solutions for large initial data is difficult. By
applying the Littlewood-Paley decomposition theory for Sobolev spaces to obtain a
losing energy estimate in $H^s$ for any $s > 2$, Wang--Xu \cite{Wang2005} obtained local solutions for any initial data and global solutions to \eqref{o}
for small initial data $(u_0, \rho_0-\bar{\rho}_0) \in H^s\times H^s$ with $s>2$. Liu--Yin \cite{Liu1,Liu2,Liu3} improved
the result of \cite{Wang2005} in the Sobolev spaces with low regularity and inhomogeneous Besov
spaces. Chen--Miao--Zhang \cite{chen0} obtained the local well-posedness of system \eqref{o} for general initial data in critical $L^p$ type Besov spaces with $1\leq p<4$ by developing a new method which relies on the smoothing properties of the heat equations. Recently, Li--Hong--Zhu \cite{lijee} proved that the system \eqref{o} is ill-posed in the critical Besov spaces with $p > 4$.

The continuous dependence is particularly important when PDEs are used to model
phenomena in the natural world since measurements are always
associated with errors. One of the first results of this type was proved by Kato \cite{Kato} who showed that the
solution operator for the (inviscid) Burgers equation is not H\"{o}lder continuous in the
$H^s(\mathbb{T})$-norm $(s > 3/2)$ for any H\"{o}lder exponent. After the phenomenon of non-uniform continuity for some dispersive equations was studied by Kenig et al. \cite{Kenig2001}, many results with regard to the non-uniform dependence on the initial data have been obtained for other nonlinear PDEs including the Euler equations \cite{HM}, the Camassa-Holm equation \cite{H-K,H-K-M,li2020}, the Benjamin-Ono equation \cite{Koch}, the compressible gas dynamics \cite{Holmes1,Keyfitz}, the Hunter-Saxton equation \cite{Holmes2} and so on. Nevertheless we notice that almost the above system mentioned is hyperbolic. As stated in \cite{Keyfitz}, the exhibition of nonuniform behavior in a hyperbolic system related to the incompressible system indicates that the nonuniform dependence is hyperbolic in nature. From the PDE's point of view, classical solution is uniform dependence on the data for a parabolic system. Naturally, one may wonder if uniform dependence on the data of solution for a coupled system which is hyperbolic-parabolic can persist or not.

In this paper, we focus on the 2D viscous shallow water equations
\begin{eqnarray}\label{vsw}
        \left\{\begin{array}{ll}
          \partial_t\rho+\D (\rho u)=0,\\
          \partial_t(\rho u)+\D(\rho u\otimes u)-\mu\D(2\rho\mathbb{D}u)+\rho\na \rho=0,\\
           \rho|_{t=0}=\rho_0,\;u|_{t=0}=u_0.
          \end{array}\right.
        \end{eqnarray}
where the strain tensor $\mathbb{D}u =\frac{1}{2}(\nabla u+\nabla^{\mathsf{T}} u)$ is the symmetric part of the velocity gradient.

For the sake of convenience, we take $\bar{\rho}_0 = 1$ and denote $\varrho=\rho-1$, we can reformulate the system \eqref{vsw} equivalently as follows
\begin{eqnarray}\label{eq}
        \left\{\begin{array}{ll}
          \partial_t\varrho+\D u+u\cd\na \varrho=-\varrho\D u,\\
          \partial_tu-{\mu} \De u-{\mu}  \na \D u+u\cd \na u+\na \varrho=2\na(\ln(1+\varrho))\cd \mathbb{D}u,\\
           \varrho|_{t=0}=\varrho_0,\;u|_{t=0}=u_0.
        \end{array}\right.
        \end{eqnarray}
Roughly speaking, the system \eqref{eq} can be regarded as a special case of compressible Navier-Stokes equations (see Refs. \cite{chen0,chen1,chen2,chen3,Danchin1,Danchin2,Danchin3,Fang,Haspot1,Haspot2,Haspot3} and the references therein).
We emphasize that the equations \eqref{eq} form a quasi-linear hyperbolic-parabolic system and possess strong nonlinear terms. Although the system \eqref{eq} is partially parabolic, owing to the first equation of \eqref{eq} which is of hyperbolic type, this creates the possibility of
the non-uniformity for the solution map. In this paper, we prove the failure of uniform dependence on the data for \eqref{eq}. To the best of our knowledge, it allows us to give a first kind of answer to the problem of the non-uniform dependence on the data for the 2D viscous shallow water equations.

Our main result of this paper is stated:
\begin{theorem}[Nonuniform dependence on initial data]\label{th1.1}
Let $s>2$.
The data-to-solution map $(u_0,\varrho_0)\mapsto \big(u(t),\varrho(t)\big)$ of the Cauchy problem \eqref{eq} is not uniformly continuous from any bounded subset in $H^{s}\times H^{s}$ into $\mathcal{C}([0,T];H^{s}\times H^{s})$, namely, there exists two sequences of solutions $\big(u_{1,n}(t),\varrho_{1,n}(t)\big)$ and $\big(u_{2,n}(t),\varrho_{2,n}(t)\big)$ such that
\begin{description}
  \item[(C.1)] $\|u_{1,n},\varrho_{1,n}\|_{H^{s}}+\|u_{2,n},\varrho_{2,n}\|_{H^{s}}\lesssim 1;$
  \item[(C.2)] $
\lim\limits_{n\to\infty}\big(\|{u}_{2,n}(0)-u_{1,n}(0)\|_{H^{s}}+\|\varrho_{2,n}(0)-\varrho_{1,n}(0)\|_{H^{s}}\big)=0;$
  \item[(C.3)] $\liminf\limits_{n\rightarrow \infty}\big(\|u_{2,n}(t)-u_{1,n}(t)\|_{H^{s}}+\|\varrho_{2,n}(t)-\varrho_{1,n}(t)\|_{H^{s}}\big)\gtrsim t$
\end{description}
for any $t\in[0,T_0]$ with small time $T_0$.
\end{theorem}
\begin{remark}\label{re1}
It should be mentioned that most non-uniformity results for solutions to other nonlinear PDEs were obtained by the Himonas-Misio{\l}ek construction in \cite{HM}. The method we used in proving Theorem \ref{th1.1} is motivated by \cite{li2020} and is completely different from that \cite{HM}.
\end{remark}
\begin{remark}\label{re2}
We note that the local well-posedness result holds for $s>1$.  However, the restriction on $s>2$ is essential in the present paper.
\end{remark}
{\bf Outline of the proof to Theorem \ref{th1.1}}\quad There are two key points in proving Theorem \ref{th1.1}:
\begin{itemize}
  \item On one hand, we construct one sequence of initial data $(u_{1,n}(0),\varrho_{1,n}(0))=(0,f_n)$, which leads to the solutions $(u_{1,n}(t),\varrho_{1,n}(t))$ to \eqref{eq} and $(u^{\rm{ap}}_{1,n},\varrho^{\rm{ap}}_{1,n})$ to the linearized system \eqref{vsw-li} respectively. Based on the special choice of $f_n$, we can prove the distance between the real solution $(u_{1,n}(t),\varrho_{1,n}(t))$ and the approximate solution $(u^{\rm{ap}}_{1,n},\varrho^{\rm{ap}}_{1,n})$ will tends to zero in $H^s$ as $n\to\infty$ . See Propositions \ref{pro1}.
  \item On the other hand, we construct another sequence of initial data $(u_{2,n}(0),\varrho_{2,n}(0))=(g_n,f_n)$, which leads to the solutions $(u_{2,n}(t),\varrho_{2,n}(t))$ to \eqref{eq} and $(u^{\rm{ap}}_{2,n},\varrho^{\rm{ap}}_{2,n})$ to the linearized system \eqref{vsw-li} respectively. Next, we aim to show that the mentioned solution $(u^{\rm{ap}}_{2,n},\varrho^{\rm{ap}}_{2,n})$ cannot approximate the real solution $(u_{2,n}(t),\varrho_{2,n}(t))$. For more details, see Propositions \ref{pro2}.
\end{itemize}
Combining the precious steps, we can conclude that the distance of two solution maps at the initial time is converging to zero, while at any later time it is bounded below by a positive constant, which means the solution maps are not uniformly continuous.\\
{\bf The structure of the paper}\quad
In Section \ref{sec2} we recall some notations and known results which will be used in the sequel. Also, we investigate the spectrum properties of the linearized system corresponding to \eqref{eq}, which will play a crucial role in the construction of approximate solutions. In Section \ref{sec3} we prove Theorem \ref{th1.1} based on the well-posedness result and some key error estimates. In Section \ref{sec4} we present some details in the computations.
\section{Preliminaries}\label{sec2}
\subsection{Notations}\label{subsec2.1}
Firstly, we introduce some notations which shall be used in this paper.
\begin{itemize}
  \item The notation $A\les B$ (resp., $A \gtrsim B$) means that there exists a harmless positive constant $c$ such that $A \leq cB$ (resp., $A \geq cB$). $A\thickapprox B$ means $A\les B$ and $A\gtrsim B$.
  \item Given a Banach space $X$, we denote its norm by $\|\cdot\|_{X}$.
  \item We shall use the simplified notation $\|f,\cdots,g\|_X=\|f\|_X+\cdots+\|g\|_X$ if there is no confusion.
  \item We shall denote by $\langle f,g\rangle$ the $L^2$ inner product of $f$ and $g$.
  \item For all $f\in \mathcal{S}'$, the Fourier transform $\mathcal{F}f$ (also denoted by $\widehat{f}$) is defined by
$$
\mathcal{F}f(\xi)=\widehat{f}(\xi)=\int_{\R^2}e^{-ix\xi}f(x)\dd x \quad\text{for any}\; \xi\in\R^2.
$$
\item We denote $\Lambda=(-\Delta)^{\fr{1}2}$. For $s\in\R$, the operator $\Lambda^s$ is defined by
$$\widehat{\Lambda^s f}(\xi)=|\xi|^{s}\widehat{f}(\xi).$$
\item The homogeneous and nonhomogeneous Sobolev space are defined by
\bbal
&\|f\|^2_{\dot{H}^s}=\|\Lambda^sf\|^2_{L^2}=\int_{\R^2}|\xi|^{2s}|\hat{f}(\xi)|^2\dd \xi,\\
&\|f\|^2_{H^s}=\int_{\R^2}(1+|\xi|^2)^s|\hat{f}(\xi)|^2\dd \xi.
\end{align*}
Then, for $s>0$, we have
$$\|f\|^2_{H^s}\thickapprox\|f\|^2_{L^2}+\|f\|^2_{\dot{H}^s}.$$
\end{itemize}

\subsection{Useful Tools}\label{subsec2.2}
Next, we review some useful tools involving the commutator and product estimates.
\begin{lemma}[See \cite{F-M-R-R}]\label{le1}
Let $s>1$ and $\nabla f,g\in H^s(\R^2)$. Then we have
\bbal
\|[\Lambda^s,f]\cdot\nabla g\|_{L^2}\leq C\|\nabla f\|_{H^s}\|g\|_{H^s}.
\end{align*}
\end{lemma}
\begin{lemma}[See \cite{Kato}]\label{le12}
Let $s>0$ and $f,g\in {\rm{Lip}}\cap H^s(\R^2)$ and $g\in L^\infty\cap H^{s-1}(\R^2)$. Then we have
\bbal
\|[\Lambda^s,f]g\|_{L^2}\leq C\big(\|\nabla f\|_{L^\infty}\|g\|_{H^{s-1}}+\|f\|_{H^s}\|g\|_{L^\infty}\big).
\end{align*}
\end{lemma}
\begin{lemma}[See \cite{B.C.D}]\label{le2}
Let $s>2$ and $f\in H^{s-2}(\R^2),g\in H^{s-1}(\R^2)$. Then we have
\bbal
\|fg\|_{H^{s-2}}\leq C\|f\|_{H^{s-2}}\|g\|_{H^{s-1}}.
\end{align*}
\end{lemma}
\begin{lemma}[See \cite{B.C.D}]\label{le3}
Let $s>0$ and $f,g\in L^\infty\cap H^s(\R^2)$. Then we have
\bbal
\|fg\|_{H^s}\leq C\big(\|f\|_{L^\infty}\|g\|_{H^s}+\|f\|_{H^s}\|g\|_{L^\infty}\big).
\end{align*}
Moreover, for $s>1$, we have the algebra estimate
\bbal
\|fg\|_{H^s}\leq C\|f\|_{H^s}\|g\|_{H^s}.
\end{align*}
\end{lemma}
\begin{lemma}[See \cite{B.C.D}]\label{le5} Let $s>0$ and $f$ be a smooth function such that $f(0)=0.$ If $u \in H^{s}\left(\mathbb{R}^{2}\right) \cap L^{\infty}\left(\mathbb{R}^{2}\right)$ then there exists a function $C$ depending only on $s$ and $f$ such that
$$
\|f(u)\|_{H^{s}} \leq C\big(\|u\|_{L^{\infty}}\big)\|u\|_{H^{s}}.
$$
\end{lemma}
\begin{lemma}[See \cite{B.C.D}]\label{le6} Let $s>1$ and $f$ be a smooth function such that $f^{\prime}(0)=0 .$ If $u, v \in H^{s}\left(\mathbb{R}^{2}\right) \cap L^{\infty}\left(\mathbb{R}^{2}\right)$, then there exists a function $C$ depending only on $s$ and $f$ such that
$$
\|f(u)-f(v)\|_{H^{s}} \leq C\big(\|u\|_{L^{\infty}},\|v\|_{L^{\infty}}\big)\|u-v\|_{H^{s}}\big(\|u\|_{H^{s}}+\|v\|_{H^{s}}\big).
$$
\end{lemma}
\subsection{The Linearized System}\label{subsec2.3}
In this subsection, borrowing the idea from \cite{chen1}, we investigate the spectrum properties of the linearized system:
\begin{eqnarray}\label{vsw-li}
        \left\{\begin{array}{ll}
          \partial_t\varrho+\D u=0,\\
          \partial_tu-{\mu} \De u-{\mu}  \na \D u+\na \varrho=0,\\
          \varrho|_{t=0}=\varrho_0,\;u|_{t=0}=u_0.
          \end{array}\right.
        \end{eqnarray}
Denote   $$d=\Lambda^{-1}\D u\quad \text{and}\quad c=\Lambda^{-1}\mathrm{curl} u,$$
then we deduce from \eqref{vsw-li}
\begin{eqnarray}\label{l1}
        \left\{\begin{array}{ll}
          \partial_t\varrho+\Lambda d=0,\\
          \partial_td-2{\mu} \Delta d-\Lambda \varrho=0,\\
          \partial_tc-{\mu} \Delta c=0
                  \end{array}\right.
        \end{eqnarray}
Let
$$\mathbf{A}=\left(
\begin{array}{cc}
0 & -|\xi| \\
|\xi| & -2\mu|\xi|^2 \\
\end{array}
\right),$$
by taking the Fourier transform of $\eqref{l1}_1$ and $\eqref{l1}_2$, we find that
\bal\label{l2}
\pa_t\left(
\begin{array}{c}
 \widehat{\varrho}\\
 \widehat{d} \\
\end{array}
\right)=\mathbf{A}\left(
\begin{array}{c}
 \widehat{\varrho}\\
 \widehat{d} \\
\end{array}
\right).
\end{align}
Straightforward calculations give the eigenvalues of the matric $\mathbf{A}$ as follows
\bbal
\lambda_{\pm}=-\mu|\xi|^2\pm \sqrt{\mu^2|\xi|^4-|\xi|^2}
\end{align*}
and
\bbal
\mathbf{T}^{-1}\mathbf{A}\mathbf{T}=\left(
\begin{array}{cc}
\lambda_{+} & 0 \\
0 & \lambda_{-} \\
\end{array}
\right),
\end{align*}
where the matrixes $\mathbf{T}$ and $\mathbf{T}^{-1}$ are given by
\bbal
\mathbf{T}=\left(
\begin{array}{cc}
\lambda_{-} & \lambda_{+} \\
-|\xi| & -|\xi| \\
\end{array}
\right)\quad\text{and}\quad \mathbf{T}^{-1}=\left(
\begin{array}{cc}
-\frac{1}{\lambda_{+}-\lambda_{-}} & -\frac{\lambda_{+}}{|\xi|(\lambda_{+}-\lambda_{-})} \\
\frac{1}{\lambda_{+}-\lambda_{-}} & \frac{\lambda_{-}}{|\xi|(\lambda_{+}-\lambda_{-})} \\
\end{array}
\right).
\end{align*}
Applying $\mathbf{T}^{-1}$ to \eqref{l2} and using Duhamel's principle, we get
\bbal
\mathbf{T}^{-1}\left(
\begin{array}{c}
 \widehat{\varrho}\\
 \widehat{d} \\
\end{array}
\right)=\left(
\begin{array}{cc}
e^{\lambda_+t}&0\\
0&e^{\lambda_-t}\\
\end{array}
\right)\mathbf{T}^{-1}\left(
\begin{array}{c}
 \widehat{\varrho}_0\\
 \widehat{d}_0 \\
\end{array}
\right),
\end{align*}
which implies that
\bal
&\widehat{\varrho}(t,\xi)=\frac{\lambda_+e^{\lambda_-t}-\lambda_-e^{\lambda_+t}}{\lambda_+-\lambda_-}\widehat{\varrho}_0
-\frac{e^{\lambda_+t}-e^{\lambda_-t}}{\lambda_+-\lambda_-}|\xi|\widehat{d}_0,\label{a}
\\&\widehat{d}(t,\xi)=
\frac{e^{\lambda_+t}-e^{\lambda_-t}}{\lambda_+-\lambda_-}|\xi|\widehat{\varrho}_0
+\frac{\lambda_+e^{\lambda_+t}-\lambda_-e^{\lambda_-t}}{\lambda_+-\lambda_-}\widehat{d}_0.\label{b}
\end{align}
From $\eqref{l1}_3$, we also have
\bal\label{c}
&\widehat{c}(t,\xi)=e^{-\mu t|\xi|^2}\widehat{c}_0.
\end{align}
Due to the vector identity
$$\Delta u=\nabla\D u+\nabla^{\bot}\mathrm{curl}u=\Lambda(\nabla d+\nabla^{\bot} c),$$
which implies
\bbal
u=-\Lambda^{-1}\na d-\Lambda^{-1}\nabla^{\bot} c,
\end{align*}
then we deduce from \eqref{b} and \eqref{c}
\bal\label{u}
\widehat{u}(t,\xi)&=
-\frac{e^{\lambda_+t}-e^{\lambda_-t}}{\lambda_+-\lambda_-}\mathcal{F}(\nabla\varrho_0)
-\frac{\lambda_+e^{\lambda_+t}-\lambda_-e^{\lambda_-t}}{\lambda_+-\lambda_-}\mathcal{F}\big(\Lambda^{-2}\na \D u_0\big)\nonumber\\
&\quad-e^{-\mu |\xi|^2t }\mathcal{F}(\Lambda^{-2}\nabla^{\bot}\mathrm{curl} u_0).
\end{align}

\section{Proof of Theorem \ref{th1.1}}\label{sec3}
This section is devoted to proving Theorem \ref{th1.1}. We begin with the well-posedness result for \eqref{vsw}.  For the details of the proof, we refer to \cite{Liu1,Liu2,Wang2005}.
\subsection{Global Well-posedness }\label{subsec3.2}
\begin{lemma}\label{le3-1}
Assume that $s>2$. For any initial data $(u_0,\varrho_0)$ which belongs to $$B_R=\big\{(\psi,\phi)\in H^s\times H^s: \|\psi\|_{H^s}+\|\phi\|_{H^s}\leq c\big\}\quad\text{for small}\;c>0,$$
then System \eqref{eq} has a unique global solution $\big(u(t),\varrho (t)\big)\in \mathcal{C}([0,+\infty);H^s\times H^s)$.\\
 Moreover, we have the solution size estimate
\begin{align}\label{s}
\|u(t)\|_{H^{\sigma}}+\|\varrho (t)\|_{H^{\sigma}}\leq C\|u_0,\varrho_0\|_{H^{\sigma}}\quad
\text{for all}\; \sigma\geq s.
\end{align}
\end{lemma}

\subsection{Approximate Solutions}\label{subsec3.3}
Let $\widehat{\phi}\in \mathcal{C}^\infty_0(\mathbb{R})$  be an even, real-valued and non-negative function on $\R$ and satisfy
\begin{numcases}{\widehat{\phi}(\xi)=}
1,&if $|\xi|\leq \frac{1}{4}$,\nonumber\\
0,&if $|\xi|\geq \frac{1}{2}$.\nonumber
\end{numcases}
Let $(u^{\rm{ap}}_{i,n},\varrho^{\rm{ap}}_{i,n})$ ($i=1,2$) be the solution of the following equation
\begin{eqnarray}\label{lin}
        \left\{\begin{array}{ll}
          \partial_t\varrho^{\rm{ap}}_{i,n}+\D u^{\rm{ap}}_{i,n}=0,\\
          \partial_tu^{\rm{ap}}_{i,n}-\mu\Delta u^{\rm{ap}}_{i,n}-\mu\na\D u^{\rm{ap}}_{i,n}+\na \varrho^{\rm{ap}}_{i,n}=0,\\
          \end{array}\right.
\end{eqnarray}
supplemented with the following initial condition, respectively,
\bal\label{c1}
(u^{\rm{ap}}_{1,n},\varrho^{\rm{ap}}_{1,n})|_{t=0}=(0,f_n)\quad\text{and}\quad (u^{\rm{ap}}_{2,n},\varrho^{\rm{ap}}_{2,n})|_{t=0}=(g_n,f_n),
\end{align}
where
\bal\label{fg}
f_n=2^{-ns}\phi(x_1)\sin (2^nx_1)\phi(x_2)\quad\text{and}\quad g_n=2^{-n}\na\big(\phi(x_1)\phi(x_2)\big),\quad n\gg1.
\end{align}
The following Lemma will play an important role in the proof of Theorem \ref{th1.1}.
\begin{lemma}\label{le4}
Let $\sigma\in\R$. Assume that $(u^{\rm{ap}}_{i,n},\varrho^{\rm{ap}}_{i,n})$ with $i=1,2$ solves System \eqref{lin}--\eqref{c1}.
Then there exists a constant $C=C(\sigma)$ such that the following statement holds
\begin{align*}
&\|\varrho^{\rm{ap}}_{1,n},u^{\rm{ap}}_{1,n}\|^2_{H^\sigma}+\mu\int^t_0\|\na u^{\rm{ap}}_{1,n},\D u^{\rm{ap}}_{1,n}\|^2_{H^\sigma}\dd \tau\leq C2^{2n(\sigma-s)},\\
&\|\varrho^{\rm{ap}}_{2,n},u^{\rm{ap}}_{2,n}\|^2_{H^\sigma}+\mu\int^t_0\|\na u^{\rm{ap}}_{2,n},\D u^{\rm{ap}}_{2,n}\|^2_{H^\sigma}\dd \tau\leq C2^{2n \max\{\sigma-s,-1\}}.
\end{align*}
\end{lemma}
{\bf Proof.}\quad Standard energy method directly gives us that
\bbal
\frac{\dd}{\dd t}\|\varrho^{\rm{ap}}_{i,n},u^{\rm{ap}}_{i,n}\|^2_{H^\sigma}+\mu\|\na u^{\rm{ap}}_{i,n}\|^2_{H^\sigma}+\mu\|\D u^{\rm{ap}}_{i,n}\|^2_{H^\sigma}=0,
\end{align*}
which implies
\bbal
&\|\varrho^{\rm{ap}}_{1,n},u^{\rm{ap}}_{1,n}\|^2_{H^\sigma}+\mu\int^t_0\|\na u^{\rm{ap}}_{1,n},\D u^{\rm{ap}}_{1,n}\|^2_{H^\sigma}\dd \tau=\|f_n\|^2_{H^\sigma},\\
&\|\varrho^{\rm{ap}}_{2,n},u^{\rm{ap}}_{2,n}\|^2_{H^\sigma}+\mu\int^t_0\|\na u^{\rm{ap}}_{2,n},\D u^{\rm{ap}}_{2,n}\|^2_{H^\sigma}\dd \tau=\|f_n,g_n\|^2_{H^\sigma}.
\end{align*}
Obviously, we have
\bal
&\|g_n\|_{H^\sigma}\les 2^{-n}\|\phi(x_1)\|_{L^2(\R)}\|\phi(x_2)\|_{L^2(\R)}\les 2^{-n},\nonumber\\
&\|f_n\|_{H^\sigma}\les 2^{-ns}\|\widetilde{f_n}(x_1,x_2)\|_{H^\sigma(\R^2)},\label{y}
\end{align}
where $\widetilde{f_n}(x_1,x_2)=\phi(x_1)\sin (2^nx_1)\phi(x_2)$.\\
Easy computations give that
\bbal
\mathcal{F}\big(\widetilde{f_n}\big)(\xi_1,\xi_2)=\frac{i}{2}\bi[\hat{\phi}(\xi_1+2^n)-\hat{\phi}(\xi_1-2^n)\bi]\hat{\phi}(\xi_2),
\end{align*}
which implies
\bbal
\mathrm{supp} \ \mathcal{F}\big(\widetilde{f_n}\big)\subset \mathcal{C}_n:=\Big\{\xi\in\R^2: \ 2^n-1\leq |\xi|\leq 2^n+1\Big\},
\end{align*}
By the definition of Sobolev space, we get
\bal\label{yyh}
\|\widetilde{f_n}(x_1,x_2)\|_{H^\sigma(\R^2)}&\les\int_{\R^2}(1+|\xi|^2)^{\sigma/2}|\mathcal{F}\big(\widetilde{f_n}\big)|^2\dd\xi\nonumber\\
&\les~2^{n\sigma}\int_{\mathcal{C}_n }|\mathcal{F}\big(\widetilde{f_n}\big)|^2\dd\xi\nonumber\\
&\les 2^{n\sigma}.
\end{align}
Inserting \eqref{yyh} into \eqref{y} enables us to finish the proof of Lemma \ref{le4}.

Before proceeding on, we give two data-to-solution maps for the Cauchy problem \eqref{eq}
\begin{align*}
&\big(u_{1,n}(0)=0,\varrho_{1,n}(0)=f_n\big)\mapsto\big(u_{1,n},\varrho_{1,n}\big),\\
&\big({u}_{2,n}(0)=g_n,{\varrho}_{2,n}(0)=f_n\big)\mapsto\big({u}_{2,n},{\varrho}_{2,n}\big),
\end{align*}
where $f_{n}$ and $g_n$ are defined in Section \ref{subsec3.3}.

Let $\ep_s=\frac12(s-2)$ and $s'=s-\ep_s>2$. For the initial data with $H^{s'}$ norm, we have
\bbal
\|u_{1,n}(0),u_{2,n}(0)\|_{H^{s'}}+\|\varrho_{1,n}(0),\varrho_{2,n}(0)\|_{H^{s'}}\leq 2^{-n\min\{\ep_s,1\}},
\end{align*}
which tends to $0$ when $n$ tends to infinity. Therefore, by Lemma \ref{le3-1}, we have the solutions size estimate which will be used implicitly in the sequel
\begin{align*}
&\|u_{1,n},\varrho_{1,n}\|_{H^\gamma}+\|u_{2,n},\varrho_{2,n}\|_{H^{\gamma}}\les2^{-n\min\{(s-\gamma),1\}}\quad
\text{for all}\; \gamma\geq s'.
\end{align*}
In particular, this implies that the condition $(\bf{C.1})$ in Theorem \ref{th1.1} holds.
\subsection{Error Estimates}\label{subsec3.4}
Letting $\varrho_{1,n}^{\rm{er}}=\varrho_{1,n}-\varrho_{1,n}^{\rm{ap}}$ and $u_{1,n}^{\rm{er}}=u_{1,n}-u_{1,n}^{\rm{ap}}$, we can find that $(\varrho_{1,n}^{\rm{er}},u_{1,n}^{\rm{er}})$ satisfies
\begin{eqnarray}\label{u1}
        \left\{\begin{array}{ll}
          \partial_t\varrho_{1,n}^{\rm{er}}+\D u_{1,n}^{\rm{er}}+u_{1,n}\cd\na \varrho_{1,n}^{\rm{er}}=\mathbf{F_1},\\
          \partial_tu_{1,n}^{\rm{er}}-{\mu} \De u_{1,n}^{\rm{er}}-{\mu}  \na \D u_{1,n}^{\rm{er}}+\na \varrho_{1,n}^{\rm{er}}+u_{1,n}\cd \na u_{1,n}^{\rm{er}}=\mathbf{F_2}+\mathbf{F_3},\\
           \varrho_{1,n}^{\rm{er}}|_{t=0}=0,\; u_{1,n}^{\rm{er}}|_{t=0}=0,       \end{array}\right.
        \end{eqnarray}
where
\bbal
&\mathbf{F_1}:=-u_{1,n}^{\rm{er}}\cd\na \varrho^{\rm{ap}}_{1,n}-\varrho_{1,n}\D u_{1,n}^{\rm{er}}-\varrho_{1,n}^{\rm{er}}\D u_{1,n}^{\rm{ap}}-\D(\varrho^{\rm{ap}}_{1,n}u_{1,n}^{\rm{ap}}),\\
&\mathbf{F_2}:=-u_{1,n}^{\rm{er}}\cd \na u^{\rm{ap}}_{1,n}- u^{\rm{ap}}_{1,n}\cd\na u^{\rm{ap}}_{1,n},\\
&\mathbf{F_3}:=2\na(\ln(1+\varrho_{1,n}))\cd \mathbb{D}u_{1,n}.
\end{align*}
The following proposition implies that the error of approximate solution and actual solution will tends to zero in $H^s$ as $n\to\infty$.
\begin{proposition}\label{pro1}
Under the assumptions of Theorem \ref{th1.1}, then we have for $t\leq 1$
\bbal
\|u_{1,n}-u_{1,n}^{\rm{ap}}\|_{H^s}+\|\varrho_{1,n}-\varrho_{1,n}^{\rm{ap}}\|_{H^s}\leq C2^{-\fr{n}{2}\min\{\ep_s,1\}}.
\end{align*}
\end{proposition}
{\bf Proof.}\quad Taking the $L^2$ inner product of $\eqref{u1}_1$ and $\eqref{u1}_2$ with $(\varrho_{1,n}^{\rm{er}},u_{1,n}^{\rm{er}})$ yields
\bal\label{y1}
&\frac12\frac{\dd }{\dd t}\|\varrho_{1,n}^{\rm{er}},u_{1,n}^{\rm{er}}\|^2_{L^2}+\mu\|\na u_{1,n}^{\rm{er}},\D u_{1,n}^{\rm{er}}\|^2_{L^{2}}\nonumber\\
=&~\fr12 \big<\D u_{1,n},|\varrho_{1,n}^{\rm{er}}|^2+|u_{1,n}^{\rm{er}}|^2\big>+ \big<\mathbf{F_1}, \varrho_{1,n}^{\rm{er}}\big>+\big<\mathbf{F_2}+\mathbf{F_3}, u_{1,n}^{\rm{er}}\big>\nonumber\\
\lesssim&~ \|u_{1,n}\|_{H^{s}}\|\varrho_{1,n}^{\rm{er}},u_{1,n}^{\rm{er}}\|^2_{L^{2}}+\|\mathbf{F_1}\|_{L^{2}}\|\varrho_{1,n}^{\rm{er}}\|_{L^2}
+\|\mathbf{F_2},\mathbf{F_3},u_{1,n}^{\rm{er}}\|^2_{L^{2}}.
\end{align}
Applying the operators $\Lambda^{s-1}\varrho_{1,n}^{\rm{er}}\Lambda^{s-1}$ and $\Lambda^{s-1}u_{1,n}^{\rm{er}}\Lambda^{s-1}$ to $\eqref{u1}_1$ and $\eqref{u1}_2$, respectively, then integrating the resulting over $\R^2$, we obtain
\bal\label{y2}
&\frac12\frac{\dd }{\dd t}\|\varrho_{1,n}^{\rm{er}},u_{1,n}^{\rm{er}}\|^2_{{\dot{H}}^{s-1}}+\mu\|\na u_{1,n}^{\rm{er}},\D u_{1,n}^{\rm{er}}\|^2_{{\dot{H}}^{s-1}}\nonumber\\
=&~\fr12 \big<\D u_{1,n},|\Lambda^{s-1}\varrho_{1,n}^{\rm{er}}|^2+|\Lambda^{s-1}u_{1,n}^{\rm{er}}|^2\big>\nonumber\\
&-\big<[\Lambda^{s-1},u_{1,n}]\cdot\nabla\varrho_{1,n}^{\rm{er}},\Lambda^{s-1}\varrho_{1,n}^{\rm{er}}\big>-\big<[\Lambda^{s-1},u_{1,n}]\cdot\nabla u_{1,n}^{\rm{er}},\Lambda^{s-1}u_{1,n}^{\rm{er}}\big>\nonumber\\
&+ \big<\Lambda^{s-1}\mathbf{F_1}, \Lambda^{s-1}\varrho_{1,n}^{\rm{er}}\big>+\big< \Lambda^{s-2}\mathbf{F_2}, \Lambda^{s}u_{1,n}^{\rm{er}}\big>+\big< \Lambda^{s-2}\mathbf{F_3}, \Lambda^{s}u_{1,n}^{\rm{er}}\big>\nonumber\\
\lesssim&~ (1+\|u_{1,n}\|_{H^{s}})\|\varrho_{1,n}^{\rm{er}},u_{1,n}^{\rm{er}}\|^2_{{{H}}^{s-1}}+\|\mathbf{F_1}\|_{{\dot{H}}^{s-1}}
\|\varrho_{1,n}^{\rm{er}}\|_{{\dot{H}}^{s-1}}\nonumber\\
&+\|\mathbf{F_2},\mathbf{F_3}\|^2_{{\dot{H}}^{s-2}}+\varepsilon\|\nabla u_{1,n}^{\rm{er}}\|^2_{{\dot{H}}^{s-1}},
\end{align}
where we have used the commutator estimate from Lemma \ref{le1}.

Combing \eqref{y1} and \eqref{y2}, we get
\bal\label{y3}
\frac{\dd }{\dd t}\|\varrho_{1,n}^{\rm{er}},u_{1,n}^{\rm{er}}\|^2_{{{H}}^{s-1}}&+\mu\|\na u_{1,n}^{\rm{er}},\D u_{1,n}^{\rm{er}}\|^2_{{{H}}^{s-1}}
\lesssim \|\varrho_{1,n}^{\rm{er}},u_{1,n}^{\rm{er}}\|^2_{{{H}}^{s-1}}\nonumber\\
&\quad+\|\mathbf{F_1}\|_{{{H}}^{s-1}}
\|\varrho_{1,n}^{\rm{er}}\|_{{{H}}^{s-1}}+\|\mathbf{F_2},\mathbf{F_3}\|^2_{{{H}}^{s-2}}+\varepsilon\|\nabla u_{1,n}^{\rm{er}}\|^2_{{\dot{H}}^{s-1}}.
\end{align}
{\bf{Estimate of $\mathbf{F_1}$.}}\quad Notice that $H^{s-1}(\R^2)$ with $s>2$ is a Banach algebra, we have
  \bbal
\|u_{1,n}^{\rm{er}}\cd\na \varrho^{\rm{ap}}_{1,n}\|_{H^{s-1}}&\les\|u_{1,n}^{\rm{er}}\|_{H^{s-1}}\|\varrho^{\rm{ap}}_{1,n}\|_{H^s}\les \|u_{1,n}^{\rm{er}}\|_{H^{s-1}},\\
\|\varrho_{1,n}\D u_{1,n}^{\rm{er}}\|_{H^{s-1}}&\les\|\varrho_{1,n}\|_{H^{s-1}}\|\D u_{1,n}^{\rm{er}}\|_{H^{s-1}}\les\|\D u_{1,n}^{\rm{er}}\|_{H^{s-1}},\\
\|\varrho_{1,n}^{\rm{er}}\D u_{1,n}^{\rm{ap}}\|_{H^{s-1}}&\les\|\varrho_{1,n}^{\rm{er}}\|_{H^{s-1}}\|u^{\rm{ap}}_{1,n}\|_{H^s}\les\|\varrho_{1,n}^{\rm{er}}\|_{H^{s-1}},\\
 \|\D(\varrho^{\rm{ap}}_{1,n}u_{1,n}^{\rm{ap}})\|_{H^{s-1}}&\les\|\varrho^{\rm{ap}}_{1,n}u_{1,n}^{\rm{ap}}\|_{H^{s}}
\les\|\varrho^{\rm{ap}}_{1,n}\|_{L^\infty}\|u_{1,n}^{\rm{ap}}\|_{H^{s}}
+\|\varrho^{\rm{ap}}_{1,n}\|_{H^{s}}\|u_{1,n}^{\rm{ap}}\|_{L^\infty}
\\&\les\|\varrho^{\rm{ap}}_{1,n}\|_{H^{s-1-\ep_s}}\|u_{1,n}^{\rm{ap}}\|_{H^{s}}
+\|\varrho^{\rm{ap}}_{1,n}\|_{H^{s}}\|u_{1,n}^{\rm{ap}}\|_{H^{s-1-\ep_s}}\les 2^{-(1+\ep_s)n},
\end{align*}
which implies
\bal\label{f1}
\|\mathbf{F_1}\|_{H^{s-1}}&\les \|u_{1,n}^{\rm{er}},\varrho_{1,n}^{\rm{er}}\|_{H^{s-1}}+\|\D u_{1,n}^{\rm{er}}\|_{H^{s-1}}+2^{-(1+\ep_s)n}.
\end{align}
{\bf{Estimate of $\mathbf{F_2}$}.}\quad Using Lemma \ref{le2} yields
\bbal
\|u_{1,n}^{\rm{er}}\cd \na u^{\rm{ap}}_{1,n}\|_{H^{s-2}}&\les\|u_{1,n}^{\rm{er}}\|_{H^{s-1}}
\|u^{\rm{ap}}_{1,n}\|_{H^{s-1}}\les\|u_{1,n}^{\rm{er}}\|_{H^{s-1}},\\
 \|u^{\rm{ap}}_{1,n}\cd\na u^{\rm{ap}}_{1,n}\|_{H^{s-2}}&\les\|u^{\rm{ap}}_{1,n}\|^2_{H^{s-1}}\les2^{-2n},
\end{align*}
which implies
\bal\label{f2}
\|\mathbf{F_2}\|_{H^{s-2}}&\les \|u_{1,n}^{\rm{er}}\|_{H^{s-1}}+2^{-2n}.
\end{align}
{\bf{Estimate of $\mathbf{F_3}$}.}\quad  We can decompose the term $\na(\ln(1+\varrho_{1,n}))\cd \mathbb{D}u_{1,n}$ as
\bbal
\na(\ln(1+\varrho_{1,n}))\cd \mathbb{D}u_{1,n}&=\na(\ln(1+\varrho_{1,n}))\cd \mathbb{D}u_{1,n}^{\rm{er}}+\na(\ln(1+\varrho_{1,n})-\ln(1+\varrho^{\rm{ap}}_{1,n}))\cd \mathbb{D}u^{\rm{ap}}_{1,n}
\\&\quad+\na(\ln(1+\varrho^{\rm{ap}}_{1,n}))\cd \mathbb{D}u^{\rm{ap}}_{1,n}\\
&=\mathbf{F_{3,1}}+\mathbf{F_{3,2}}+\mathbf{F_{3,3}}.
\end{align*}
For the first two terms, by Lemma \ref{le2} and Lemmas \ref{le5}--\ref{le6}, we have
\bbal
&\|\mathbf{F_{3,1}}\|_{H^{s-2}}\les\|\mathbb{D}u_{1,n}^{\rm{er}}\|_{H^{s-2}}\|\varrho_{1,n}\|_{H^s}\les\|u_{1,n}^{\rm{er}}\|_{H^{s-1}},
\\&\|\mathbf{F_{3,2}}\|_{H^{s-2}}\les\|\varrho_{1,n}^{\rm{er}}\|_{H^{s-1}}\|u^{\rm{ap}}_{1,n}\|_{H^s}\les\|\varrho_{1,n}^{\rm{er}}\|_{H^{s-1}}.
\end{align*}
For the third term, by Lemmas \ref{le3} and \ref{le5}, we have
\bbal
 \|\mathbf{F_{3,3}}\|_{H^{s-2}}&\les\|\na(\ln(1+\varrho^{\rm{ap}}_{1,n}))\|_{L^\infty}\|u_{1,n}^{\rm{ap}}\|_{H^{s-1}}
+\|\na(\ln(1+\varrho^{\rm{ap}}_{1,n}))\|_{H^{s-2}}\|\na u_{1,n}^{\rm{ap}}\|_{L^\infty}
\\&\les\|\varrho^{\rm{ap}}_{1,n}\|_{H^{s-\ep_s}}\|u_{1,n}^{\rm{ap}}\|_{H^{s-1}}
+\|\varrho^{\rm{ap}}_{1,n}\|_{H^{s-1}}\|u_{1,n}^{\rm{ap}}\|_{H^{s-\ep_s}}\les 2^{-(1+\ep_s)n},
\end{align*}
which implies
\bal\label{f3}
\|\mathbf{F_3}\|_{H^{s-2}}&\les \|u_{1,n}^{\rm{er}},\varrho_{1,n}^{\rm{er}}\|_{H^{s-1}}+2^{-(1+\ep_s)n}.
\end{align}
Putting the above estimates \eqref{f1}--\eqref{f3} together with \eqref{y3}, we can obtain
\bbal
\frac{\dd }{\dd t}\|\varrho_{1,n}^{\rm{er}},u_{1,n}^{\rm{er}}\|^2_{{{H}}^{s-1}}+\mu\|\na u_{1,n}^{\rm{er}},\D u_{1,n}^{\rm{er}}\|^2_{{{H}}^{s-1}}
\les&~ \|\varrho_{1,n}^{\rm{er}},u_{1,n}^{\rm{er}}\|^2_{H^{s-1}}+2^{-2n(1+\min\{\ep_s,1\})}\\&+\varepsilon\|\na u_{1,n}^{\rm{er}},\D u_{1,n}^{\rm{er}}\|^2_{{{H}}^{s-1}}.
\end{align*}
Absorbing the $\varepsilon$-term and using Gronwall's inequality yield
\bbal
\|\varrho_{1,n}^{\rm{er}},u_{1,n}^{\rm{er}}\|_{{{H}}^{s-1}}
\leq C2^{-n(1+\min\{\ep_s,1\})}.
\end{align*}
An interpolation argument leads to
\bbal
\|\varrho_{1,n}^{\rm{er}},u_{1,n}^{\rm{er}}\|_{{{H}}^{s}}\lesssim
\|\varrho_{1,n}^{\rm{er}},u_{1,n}^{\rm{er}}\|_{{{H}}^{s-1}}^{\frac{1}{2}}\|\varrho_{1,n}^{\rm{er}},u_{1,n}^{\rm{er}}\|_{{{H}}^{s+1}}^{\frac{1}{2}}\lesssim 2^{-\fr{n}{2}\min\{\ep_s,1\}}.
\end{align*}
Thus, we have finished the proof of Proposition \ref{pro1}.

Denoting $V^{\rm{ap}}_{n}=-u_{2,n}^{\rm{ap}}\cd\na \varrho_{2,n}^{\rm{ap}}$ and introducing the errors $$\varrho_{2,n}^{\rm{er}}=\varrho_{2,n}-\varrho_{2,n}^{\rm{ap}}-tV^{\rm{ap}}_{n}\quad\text{and}\quad u_{2,n}^{\rm{er}}=u_{2,n}-u_{2,n}^{\rm{ap}},$$ we can find that $(\varrho_{2,n}^{\rm{er}},u_{2,n}^{\rm{er}})$ satisfies
\begin{eqnarray}\label{CNS-re}
        \left\{\begin{array}{ll}
          \partial_t\varrho_{2,n}^{\rm{er}}+\D u_{2,n}^{\rm{er}}+u_{2,n}\cd\na \varrho_{2,n}^{\rm{er}}=\mathbf{G_1}
        -t\mathbf{G_4}-t\pa_tV^{\rm{ap}}_{n},\\
          \partial_tu_{2,n}^{\rm{er}}-{\mu} \De u_{2,n}^{\rm{er}}-{\mu}  \na \D u_{2,n}^{\rm{er}}+\na \varrho_{2,n}^{\rm{er}}+u_{2,n}\cd \na u_{2,n}^{\rm{er}}=\mathbf{G_2}+\mathbf{G_3}-t\na V^{\rm{ap}}_{n}, \\
           \varrho_{2,n}^{\rm{er}}|_{t=0}=0,\; u_{2,n}^{\rm{er}}|_{t=0}=0,         \end{array}\right.
        \end{eqnarray}
where
\bbal
&\mathbf{G_1}:=-u_{2,n}^{\rm{er}}\cd\na \varrho^{\rm{ap}}_{2,n}-\varrho_{2,n}\D u_{2,n}^{\rm{er}}-\varrho_{2,n}^{\rm{er}}\D u_{2,n}^{\rm{ap}}-\varrho^{\rm{ap}}_{2,n}\D u_{2,n}^{\rm{ap}},\\
&\mathbf{G_2}:=-u_{2,n}^{\rm{er}}\cd \na u^{\rm{ap}}_{2,n}- u^{\rm{ap}}_{2,n}\cd\na u^{\rm{ap}}_{2,n},\\
&\mathbf{G_3}:=2\na(\ln(1+\varrho_{2,n}))\cd \mathbb{D}u_{2,n},\\
&\mathbf{G_4}:=u_{2,n}^{\rm{er}}\cd\na V^{\rm{ap}}_{n}+\D(V^{\rm{ap}}_{n} u_{2,n}^{\rm{ap}}).
\end{align*}

\begin{proposition}\label{pro2}
Under the assumptions of Theorem \ref{th1.1}, then we have  for $t\leq 1$
\bbal
\|u_{2,n}-u_{2,n}^{\rm{ap}}\|_{H^s}+\|\varrho_{2,n}-\varrho_{2,n}^{\rm{ap}}-tV^{\rm{ap}}_{n}\|_{H^s}\leq C2^{-n\min\{\ep_s,1\}}+Ct^{\fr32}.
\end{align*}
\end{proposition}
{\bf Proof.}\quad  Following the procedure in \eqref{y1} and \eqref{y2}, we obtain
\bal\label{y4}
&\frac{\dd }{\dd t}\|\varrho_{2,n}^{\rm{er}},u_{2,n}^{\rm{er}}\|^2_{{{H}}^{s}}+\mu\|\na u_{2,n}^{\rm{er}},\D u_{2,n}^{\rm{er}}\|^2_{{{H}}^{s}}\nonumber\\
\lesssim&~ (1+\|u_{2,n}\|_{H^{s}})\|\varrho_{2,n}^{\rm{er}},u_{2,n}^{\rm{er}}\|^2_{{{H}}^{s}}\nonumber\\
&+\big(\|\mathbf{G_1}\|_{{{H}}^{s}}+\|t\pa_tV^{\rm{ap}}_{n}\|_{{{H}}^{s}}+t\|\mathbf{G_4}\|_{{{H}}^{s}}\big)
\|\varrho_{2,n}^{\rm{er}}\|_{{{H}}^{s}}\nonumber\\
&+\|\mathbf{G_2},\mathbf{G_3},t\nabla V^{\rm{ap}}_{n}\|^2_{{{H}}^{s-1}}\nonumber\\
\lesssim&~ \|\varrho_{2,n}^{\rm{er}},u_{2,n}^{\rm{er}}\|^2_{{{H}}^{s}}
+\|\mathbf{G_1}\|_{{{H}}^{s}}\|\varrho_{2,n}^{\rm{er}}\|_{{{H}}^{s}}+\|\mathbf{G_2},\mathbf{G_3}\|^2_{{{H}}^{s-1}}+\|t\pa_tV^{\rm{ap}}_{n}\|^2_{{{H}}^{s}}
\nonumber\\
&+t^2\big(\|\mathbf{G_4}\|^2_{{{H}}^{s}}+\|V^{\rm{ap}}_{n}\|^2_{{{H}}^{s}}\big).
\end{align}
We should mentioned that the commutator estimate from Lemma \ref{le12} was used here.

It is not hard to deduce that for $k\in\{-1,0,1\}$
\bbal
\|V^{\rm{ap}}_{n}\|_{H^{s+k}}&\les \|u_{2,n}^{\rm{ap}}\|_{L^\infty}\|\na \varrho_{2,n}^{\rm{ap}}\|_{H^{s+k}}
+\|u_{2,n}^{\rm{ap}}\|_{H^{s+k}}\|\na \varrho_{2,n}^{\rm{ap}}\|_{L^\infty}
\\&\les \|u_{2,n}^{\rm{ap}}\|_{H^{s-1}}\|\na \varrho_{2,n}^{\rm{ap}}\|_{H^{s+k}}
+\|u_{2,n}^{\rm{ap}}\|_{H^{s+k}}\|\na \varrho_{2,n}^{\rm{ap}}\|_{H^{s-1}}\\&\les 2^{kn}.
\end{align*}
{\bf{Estimate of $\mathbf{G_1}$.}}\quad Using Lemma \ref{le3} yields
\bbal
\|u_{2,n}^{\rm{er}}\cd\na \varrho^{\rm{ap}}_{2,n}\|_{H^{s}}&\les\|u_{2,n}^{\rm{er}}\|_{L^\infty}\|\varrho^{\rm{ap}}_{2,n}\|_{H^{s+1}}+\|u_{2,n}^{\rm{er}}\|_{H^s}
\|\nabla\varrho^{\rm{ap}}_{2,n}\|_{L^\infty}\\
&\les 2^n\|u_{2,n}^{\rm{er}}\|_{H^{s-1}}+\|u_{2,n}^{\rm{er}}\|_{H^{s}},\\
\|\varrho_{2,n}\D u_{2,n}^{\rm{er}}\|_{H^{s}}&\les\|\D u_{2,n}^{\rm{er}}\|_{H^{s}},\\
\|\varrho_{2,n}^{\rm{er}}\D u_{2,n}^{\rm{ap}}\|_{H^{s}}&\les 2^n\|\varrho_{2,n}^{\rm{er}}\|_{H^{s-1}}+\|\varrho_{2,n}^{\rm{er}}\|_{H^{s}},\\
 \|\varrho^{\rm{ap}}_{2,n}\D u_{2,n}^{\rm{ap}}\|_{H^{s}}&\les\|\varrho_{2,n}^{\rm{ap}}\|_{H^{s-1}}\|\D u^{\rm{ap}}_{2,n}\|_{H^{s}}+\|\varrho_{2,n}^{\rm{ap}}\|_{H^s}
\|u^{\rm{ap}}_{2,n}\|_{H^{s-\varepsilon_s}}\\
&\les 2^{-n}\|\D u^{\rm{ap}}_{2,n}\|_{H^{s}}+2^{-n\min\{\ep_s,1\}},
\end{align*}
which implies that
\bal\label{g1}
 \|\mathbf{G_{1}}\|_{H^{s}}&\les2^n\|u_{2,n}^{\rm{er}},\varrho_{2,n}^{\rm{er}}\|_{H^{s-1}}+\|\D u_{2,n}^{\rm{er}}\|_{H^{s}}+\|u_{2,n}^{\rm{er}},\varrho_{2,n}^{\rm{er}}\|_{H^{s}}\nonumber\\
 &\quad+2^{-n}\|\D u^{\rm{ap}}_{2,n}\|_{H^{s}}+2^{-n\min\{\ep_s,1\}}.
\end{align}
{\bf{Estimate of $\mathbf{G_2}$}.}\quad Notice that $H^{s-1}(\R^2)$ with $s>2$ is a Banach algebra, we have
\bbal
\|u_{2,n}^{\rm{er}}\cd \na u^{\rm{ap}}_{2,n}\|_{H^{s-1}}&\les \|u_{2,n}^{\rm{er}}\|_{H^{s-1}}\|u_{2,n}^{\rm{ap}}\|_{H^{s}}\les \|u_{2,n}^{\rm{er}}\|_{H^{s}},\\
\|u^{\rm{ap}}_{2,n}\cd\na u^{\rm{ap}}_{2,n}\|_{H^{s-1}}&\les\|u^{\rm{ap}}_{2,n}\|_{H^{s-1}}\|u_{2,n}^{\rm{ap}}\|_{H^{s}}
\les 2^{-n},
\end{align*}
which implies that
\bal\label{g2}
 \|\mathbf{G_{2}}\|_{H^{s-1}}&\les\|u_{2,n}^{\rm{er}}\|_{H^{s}}+2^{-n}.
\end{align}
{\bf{Estimate of $\mathbf{G_3}$}.}\quad
For the term $\na(\ln(1+\varrho_{2,n}))\cd \mathbb{D}u_{2,n}$, we can decompose it as
\bbal
\na(\ln(1+\varrho_{2,n}))\cd \mathbb{D}u_{2,n}&=\na(\ln(1+\varrho_{2,n}))\cd \mathbb{D}u_{2,n}^{\rm{er}}\\&\quad
+\na(\ln(1+\varrho_{2,n})-\ln(1+\varrho^{\rm{ap}}_{2,n}+tV^{\rm{ap}}_{n}))\cd \mathbb{D}u^{\rm{ap}}_{2,n}
\\&\quad+\na(\ln(1+\varrho^{\rm{ap}}_{2,n}+tV^{\rm{ap}}_{n}))\cd \mathbb{D}u^{\rm{ap}}_{2,n}\\&
=\mathbf{G_{3,1}}+\mathbf{G_{3,2}}+\mathbf{G_{3,3}}.
\end{align*}
By Lemmas \ref{le2} and \ref{le3}, it is easy to deduce that
\bbal
&\|\mathbf{G_{3,1}}\|_{H^{s-1}}\lesssim \|u_{2,n}^{\rm{er}}\|_{H^{s}},\\
&\|\mathbf{G_{3,2}}\|_{H^{s-1}}\lesssim \|\varrho_{2,n}^{\rm{er}}\|_{H^{s}},\\
&\|\mathbf{G_{3,3}}\|_{H^{s-1}}\les2^{-n\min\{\ep_s,1\}},
\end{align*}
which gives that
\bal\label{g3}
 \|\mathbf{G_{3}}\|_{H^{s-1}}&\les\|u_{2,n}^{\rm{er}},\varrho_{2,n}^{\rm{er}}\|_{H^{s}}+2^{-n\min\{\ep_s,1\}}.
\end{align}
{\bf{Estimate of $\mathbf{G_4}$}.}\quad By Lemma \ref{le3} again, we have
\bbal
\|u_{2,n}^{\rm{er}}\cd\na V^{\rm{ap}}_{n}\|_{H^{s}}&\les\|u_{2,n}^{\rm{er}}\|_{L^\infty}\|V^{\rm{ap}}_{n}\|_{H^{s+1}}+\|u_{2,n}^{\rm{er}}\|_{H^s}
\|\nabla V^{\rm{ap}}_{n}\|_{L^\infty}\les 2^n\|u_{2,n}^{\rm{er}}\|_{H^{s-1}}+\|u_{2,n}^{\rm{er}}\|_{H^{s}},\\
 \|\D(V^{\rm{ap}}_{2,n}u_{2,n}^{\rm{ap}})\|_{H^{s}}&\les \|V^{\rm{ap}}_{2,n}\D u_{2,n}^{\rm{ap}}\|_{H^{s}}+ \|u_{2,n}^{\rm{ap}}\cdot\nabla V^{\rm{ap}}_{2,n}\|_{H^{s}}\les 1
\end{align*}
and
\bbal
\|u_{2,n}^{\rm{er}}\cd\na V^{\rm{ap}}_{n}\|_{H^{s-1}}&\les\|u_{2,n}^{\rm{er}}\|_{L^\infty}\|V^{\rm{ap}}_{n}\|_{H^{s}}+\|u_{2,n}^{\rm{er}}\|_{H^{s-1}}
\|\nabla V^{\rm{ap}}_{n}\|_{L^\infty}\les \|u_{2,n}^{\rm{er}}\|_{H^{s-1}},\\
 \|\D(V^{\rm{ap}}_{2,n}u_{2,n}^{\rm{ap}})\|_{H^{s-1}}&\les \|V^{\rm{ap}}_{2,n}u_{2,n}^{\rm{ap}}\|_{H^{s}}\les 2^{-n},
\end{align*}
from which, we obtain
 \bal\label{g4}
 \|\mathbf{G_{4}}\|_{H^{s}}&\les2^n\|u_{2,n}^{\rm{er}}\|_{H^{s-1}}+\|u_{2,n}^{\rm{er}}\|_{H^{s}}+1\quad\text{and}\quad\|\mathbf{G_{4}}\|_{H^{s-1}}\les\|u_{2,n}^{\rm{er}}\|_{H^{s-1}}+2^{-n}.
\end{align}
Next, we need to deal with the involved term $\|t\pa_tV^{\rm{ap}}_{n}\|_{{{H}}^{s}}$. \\
{\bf{Estimate of $\|t\pa_tV^{\rm{ap}}_{n}\|_{{{H}}^{s}}$}.}\quad Direct calculation gives that
\begin{align*}
&-\pa_tV_{n}^{\rm{ap}}=\pa_tu_{2,n}^{\rm{ap}}\cdot\nabla\varrho_{2,n}^{\rm{ap}}+u_{2,n}^{\rm{ap}}\cdot\nabla\pa_t\varrho_{2,n}^{\rm{ap}},
\end{align*}
using the following estimates whose proof are relegated to {\bf{A.1--A.2}} in Section \ref{sec4}
\bal
&\|t\pa_tu_{2,n}^{\rm{ap}}\|_{H^{s+k}}\leq C2^{(k-1)n}+Ct2^{-n},\label{a1}
\\&\|\pa_t\varrho_{2,n}^{\rm{ap}}\|_{H^{s+k}}\leq C2^{kn},\label{a2}
\end{align}
then we can estimate $\|t\pa_tV^{\rm{ap}}_{n}\|_{{{H}}^{s}}$ as
\bbal
\|t\pa_tV^{\rm{ap}}_{n}\|_{H^{s+k}}&\les \|t\pa_tu_{2,n}^{\rm{ap}}\|_{L^\infty}\|\na \varrho_{2,n}^{\rm{ap}}\|_{H^{s+k}}
+\|t\pa_tu_{2,n}^{\rm{ap}}\|_{H^{s+k}}\|\na \varrho_{2,n}^{\rm{ap}}\|_{L^\infty}
\\&\quad +t\|u_{2,n}^{\rm{ap}}\|_{L^\infty}\|\na \pa_t\varrho_{2,n}^{\rm{ap}}\|_{H^{s+k}}
+t\|u_{2,n}^{\rm{ap}}\|_{H^{s+k}}\|\na \pa_t\varrho_{2,n}^{\rm{ap}}\|_{L^\infty}
\\&\les \|t\pa_tu_{2,n}^{\rm{ap}}\|_{H^{s-1}}\|\varrho_{2,n}^{\rm{ap}}\|_{H^{s+k+1}}
+\|t\pa_tu_{2,n}^{\rm{ap}}\|_{H^{s+k}}\|\varrho_{2,n}^{\rm{ap}}\|_{H^{s}}
\\&\quad +t\|u_{2,n}^{\rm{ap}}\|_{H^{s-1}}\|\pa_t\varrho_{2,n}^{\rm{ap}}\|_{H^{s+k+1}}
+t\|u_{2,n}^{\rm{ap}}\|_{H^{s+k}}\|\pa_t\varrho_{2,n}^{\rm{ap}}\|_{H^{s}}
\\&\les t2^{kn}+2^{(k-1)n}.
\end{align*}
Putting the above estimates into \eqref{y3} yields
\bal\label{a3}
\frac{\dd }{\dd t}\|\varrho_{2,n}^{\rm{er}},u_{2,n}^{\rm{er}}\|^2_{{{H}}^{s}}
\les&~\|\varrho_{2,n}^{\rm{er}},u_{2,n}^{\rm{er}}\|^2_{{{H}}^{s}}+2^{2n}\|\varrho_{2,n}^{\rm{er}},u_{2,n}^{\rm{er}}\|^2_{H^{s-1}}+2^{-2n}\|\D u_{2,n}^{\rm{ap}}\|^2_{H^{s}}\nonumber\\
&+2^{-2n\min\{\ep_s,1\}}+t^2.
\end{align}
Finally, to close the above, we have to estimate $\|\varrho_{2,n}^{\rm{er}},u_{2,n}^{\rm{er}}\|^2_{H^{s-1}}$. Following the procedure in \eqref{y1} and \eqref{y2} once again, we obtain
\bal\label{y5}
&\frac{\dd }{\dd t}\|\varrho_{2,n}^{\rm{er}},u_{2,n}^{\rm{er}}\|^2_{{{H}}^{s-1}}+\mu\|\na u_{2,n}^{\rm{er}},\D u_{2,n}^{\rm{er}}\|^2_{{{H}}^{s-1}}\nonumber\\
\lesssim&~ (1+\|u_{2,n}\|_{H^{s}})\|\varrho_{2,n}^{\rm{er}},u_{2,n}^{\rm{er}}\|^2_{{{H}}^{s-1}}\nonumber\\
&+\big(\|\mathbf{G_1}\|_{{{H}}^{s-1}}+\|t\pa_tV^{\rm{ap}}_{n}\|_{{{H}}^{s-1}}\big)\|\varrho_{2,n}^{\rm{er}}\|_{H^{s-1}}+\|\mathbf{G_2},\mathbf{G_3}\|^2_{{{H}}^{s-2}}
\nonumber\\
&+t^2\big(\|\mathbf{G_4}\|^2_{H^{s-1}}+\|V^{\rm{ap}}_{n}\|^2_{{{H}}^{s-1}}\big)\nonumber\\
\lesssim&~\|\varrho_{2,n}^{\rm{er}},u_{2,n}^{\rm{er}}\|^2_{{{H}}^{s-1}}+
\|\mathbf{G_1}\|_{{{H}}^{s-1}}\|\varrho_{2,n}^{\rm{er}}\|_{H^{s-1}}+\|\mathbf{G_2},\mathbf{G_3}\|^2_{{{H}}^{s-2}}+\|t\pa_tV^{\rm{ap}}_{n}\|^2_{{{H}}^{s-1}}
\nonumber\\
&+t^2\big(\|\mathbf{G_4}\|^2_{H^{s-1}}+\|V^{\rm{ap}}_{n}\|^2_{{{H}}^{s-1}}\big).
\end{align}
It should be mentioned that the first three terms of $\mathbf{G_1}$ can be done as that of $\mathbf{F_1}$. The only difference is the forth term of $\mathbf{G_1}$. In fact, we estimate it as
\bbal
\|\varrho^{\rm{ap}}_{2,n}\D u_{2,n}^{\rm{ap}}\|_{{{H}}^{s-1}}\les\|\varrho^{\rm{ap}}_{2,n}\|_{{{H}}^{s-1}}\|\D u_{2,n}^{\rm{ap}}\|_{{{H}}^{s-1}}\les2^{-n}\|\D u_{2,n}^{\rm{ap}}\|_{{{H}}^{s-1}}.
\end{align*}
Then we have
\bbal
\|\mathbf{G_1}\|_{H^{s-1}}&\les \|u_{1,n}^{\rm{er}},\varrho_{1,n}^{\rm{er}}\|_{H^{s-1}}+\|\D u_{1,n}^{\rm{er}}\|_{H^{s-1}}+2^{-n}\|\D u_{2,n}^{\rm{ap}}\|_{{{H}}^{s-1}}.
\end{align*}
The estimates of $\mathbf{G_2}$ and $\mathbf{G_3}$ in $H^{s-2}$ can be done as that of $\mathbf{F_2}$ and $\mathbf{F_3}$, respectively. We have
  \bbal
\|\mathbf{G_2},\mathbf{G_3}\|_{H^{s-2}}&\les \|u_{1,n}^{\rm{er}},\varrho_{1,n}^{\rm{er}}\|_{H^{s-1}}+2^{-(1+\ep_s)n}+2^{-2n}.
\end{align*}
Combining the above estimates with \eqref{y5}, we can obtain
\bbal
&\frac{\dd }{\dd t}\|\varrho_{2,n}^{\rm{er}},u_{2,n}^{\rm{er}}\|^2_{{{H}}^{s-1}}
\lesssim \|\varrho_{2,n}^{\rm{er}},u_{2,n}^{\rm{er}}\|^2_{{{H}}^{s-1}}+2^{-2n}\|\D u_{2,n}^{\rm{ap}}\|^2_{{{H}}^{s-1}}+2^{-2n(1+\min\{\ep_s,1\})}+2^{-2n}t^2,
\end{align*}
which follows from Gronwall's inequality and Lemma \ref{le4} that
\bal\label{er}
\|\varrho_{2,n}^{\rm{er}},u_{2,n}^{\rm{er}}\|^2_{{{H}}^{s-1}}
\les2^{-2n(1+\min\{\ep_s,1\})}+2^{-2n}t^3.
\end{align}
Plugging \eqref{er} into \eqref{a3} yields
\bbal
\frac{\dd }{\dd t}\|\varrho_{2,n}^{\rm{er}},u_{2,n}^{\rm{er}}\|^2_{{{H}}^{s}}
\les&~\|\varrho_{2,n}^{\rm{er}},u_{2,n}^{\rm{er}}\|^2_{{{H}}^{s}}+2^{-2n}\|\D u_{2,n}^{\rm{ap}}\|^2_{H^{s}}+2^{-2n\min\{\ep_s,1\}}+t^2.
\end{align*}
Using Gronwall's inequality and Lemma \ref{le4} yields
\bbal
\|\varrho_{2,n}^{\rm{er}},u_{2,n}^{\rm{er}}\|^2_{{{H}}^{s}}
\les&~2^{-2n\min\{\ep_s,1\}}+t^3.
\end{align*}
Thus, we have finished the proof of Proposition \ref{pro2}.
\subsection{Non-uniform Continuous Dependence}\label{subsec3.5}
With Propositions \ref{pro1}--\ref{pro2} in hand, we can prove Theorem \ref{th1.1}.

{\bf Behavior at time $t=0$.}\quad
Obviously, we have
\begin{align*}
\|{u}_{2,n}(0)-u_{1,n}(0)\|_{H^{s}}+\|\varrho_{2,n}(0)-\varrho_{1,n}(0)\|_{H^{s}}=\|g_n\|_{H^{s}}\leq C2^{-n},
\end{align*}
which means that the condition $(\bf{C.2})$ in Theorem \ref{th1.1} holds.

{\bf Behavior at time $t>0$.}\quad
Notice that
\begin{align*}
&u_{2,n}-u_{1,n}=u^{\rm{er}}_{2,n}-u^{\rm{er}}_{1,n}+u^{\rm{ap}}_{2,n}-u^{\rm{ap}}_{1,n},\nonumber\\
&\varrho_{2,n}-\varrho_{1,n}=t{V}_n^{\rm{ap}}+\varrho^{\rm{er}}_{2,n}-\varrho^{\rm{er}}_{1,n}+\varrho^{\rm{ap}}_{2,n}-\varrho^{\rm{ap}}_{1,n},
\end{align*}
by the triangle inequality and Propositions \ref{pro1} and \ref{pro2}, we deduce that  for $t\leq 1$
\begin{align}\label{z}
&\|\varrho_{2,n}(t)-\varrho_{1,n}(t)\|_{H^{s}}+\|u_{2,n}(t)-u_{1,n}(t)\|_{H^{s}}\nonumber\\
\geq&~ t\|V^{\rm{ap}}_{n}\|_{H^s}-\big(\|u^{\rm{ap}}_{2,n}-u^{\rm{ap}}_{1,n}\|_{H^s}+\|\varrho^{\rm{ap}}_{2,n}-\varrho^{\rm{ap}}_{1,n}\|_{H^s}\big)
-\big(\|u^{\rm{er}}_{1,n},\varrho^{\rm{er}}_{1,n}\|_{H^s}+\|u^{\rm{er}}_{2,n},\varrho^{\rm{er}}_{2,n}\|_{H^s}\big)\nonumber\\
\gtrsim&~  t\|V^{\rm{ap}}_{n}\|_{H^s}-2^{-n}-2^{-\frac{n}{2}\min\{\ep_s,1\}}-2^{-n\min\{\ep_s,1\}}-t^{\fr32},
\end{align}
where we have used
$$\|u^{\rm{ap}}_{2,n}-u^{\rm{ap}}_{1,n}\|^2_{H^s}+\|\varrho^{\rm{ap}}_{2,n}-\varrho^{\rm{ap}}_{1,n}\|^2_{H^s}\les\|g_n\|^2_{H^s}\les 2^{-2n}.$$
Notice that $-V^{\rm{ap}}_{n}=u_{2,n}^{\rm{ap}}\cd\na \varrho_{2,n}^{\rm{ap}}$, we decompose it as
\begin{align*}
-V^{\rm{ap}}_{n}=(u_{2,n}^{\rm{ap}})_1\pa_1 \varrho_{2,n}^{\rm{ap}}+(u_{2,n}^{\rm{ap}})_2\pa_2 \varrho_{2,n}^{\rm{ap}}=(g_n)_1\pa_1f_n+E,
\end{align*}
where
$$E:=(u_{2,n}^{\rm{ap}}-g_n)_1\pa_1\varrho_{2,n}^{\rm{ap}}+(g_n)_1\pa_1(\varrho_{2,n}^{\rm{ap}}-f_n)+(u_{2,n}^{\rm{ap}})_2\pa_2\varrho_{2,n}^{\rm{ap}}.$$
Then we have  for $t\leq 1$
\bal\label{E}
\|E\|_{H^{s}}\leq Ct+C2^{-n}.
\end{align}
For more details of proof, see {\bf{A.3}} in Section \ref{sec4}.

Furthermore, \eqref{z} reduces to
\begin{align}\label{z1}
\|\varrho_{2,n}(t)-\varrho_{1,n}(t)\|_{H^{s}}+\|u_{2,n}(t)-u_{1,n}(t)\|_{H^{s}}
\gtrsim t\|(g_n)_1\pa_1f_n\|_{H^s}-2^{-\frac{n}{2}\min\{\ep_s,1\}}-t^{\fr32},
\end{align}
combining the following estimate whose proof is postponed to {\bf{A.4}} in Section \ref{sec4}
\begin{align}\label{zw}
\liminf_{n\rightarrow \infty}\|(g_n)_1\pa_1f_n\|_{H^s}\gtrsim \|\phi'\phi\|^2_{L^2}\|\phi^2\|_{L^2},
\end{align}
we get from \eqref{z1} that
\bbal
\liminf_{n\rightarrow \infty}\big(\|\varrho_{2,n}(t)-\varrho_{1,n}(t)\|_{H^{s}}+\|u_{2,n}(t)-u_{1,n}(t)\|_{H^{s}}\big)\gtrsim t\quad\text{for} \ t \ \text{small enough},
\end{align*}
which is nothing but the condition $(\bf{C.3})$ in Theorem \ref{th1.1}.

Thus, we complete the proof of Theorem \ref{th1.1}.
\section{Appendix}\label{sec4}
\setcounter{equation}{0}
For the sake of convenience, here we present more details in the computations.\\
{\bf A.1\quad Proof of \eqref{a1}.}\quad
When $|\xi|\rightarrow\infty$, we have
\bal\label{lambda-es}
\lambda_-(\xi)\sim -2\mu|\xi|^2 \quad \text{and}\quad\lambda_+(\xi)\sim -\frac{1}{2\mu}.
\end{align}
From \eqref{u}, then
\bbal
\mathcal{F}\big(u_{2,n}^{\rm{ap}}\big)(t,\xi)&=
-\frac{e^{\lambda_+t}-e^{\lambda_-t}}{\lambda_+-\lambda_-}\mathcal{F}\big(\nabla f_n\big)
-\frac{\lambda_+e^{\lambda_+t}-\lambda_-e^{\lambda_-t}}{\lambda_+-\lambda_-}\mathcal{F}\big(\Lambda^{-2}\na \D g_n\big)
\end{align*}
which gives
\bbal
&\mathcal{F}\big(t\pa_tu_{2,n}^{\rm{ap}}\big)=
-\frac{t\lambda_+e^{\lambda_+t}-t\lambda_-e^{\lambda_-t}}{\lambda_+-\lambda_-}\mathcal{F}\big(\nabla f_n\big)
-t\frac{\lambda_{+}^2e^{\lambda_+t}-\lambda_{-}^2e^{\lambda_-t}}{\lambda_+-\lambda_-}\mathcal{F}\big(\Lambda^{-2}\na \D g_n\big).
\end{align*}
Using the simple facts
\bbal
&|t\lambda_{\pm}e^{\lambda_{\pm}t}|\leq1\quad\text{if }\;\xi\in\mathcal{C}_n \\
&\Big|\frac{\lambda_{+}^2e^{\lambda_+t}-\lambda_{-}^2e^{\lambda_-t}}{\lambda_+-\lambda_-}
\Big|\lesssim |\lambda_++\lambda_-|+t|\lambda_-|^2\les 1\quad\text{if }\;\xi\in\mathcal{B}(0,1),
\end{align*}
then we have
\bbal
\big|\mathcal{F}\big(t\pa_tu_{2,n}^{\rm{ap}}\big)\big|&\leq
\Big|\frac{t\lambda_+e^{\lambda_+t}-t\lambda_-e^{\lambda_-t}}{\lambda_+-\lambda_-}\mathcal{F}\big(\nabla f_n\big)\Big|
+t\Big|\frac{\lambda_{+}^2e^{\lambda_+t}-\lambda_{-}^2e^{\lambda_-t}}{\lambda_+-\lambda_-}\mathcal{F}\big(\Lambda^{-2}\na \D g_n\big)\Big|\\
&\leq
\frac{2}{|\lambda_+-\lambda_-|}\big|\mathcal{F}\big(\nabla f_n\big)\big|
+t\big|\mathcal{F}\big(\Lambda^{-2}\na \D g_n\big)\big|\\
&\approx\big|\mathcal{F}\big(\Lambda^{-2}\nabla f_n\big)\big|
+t\big|\mathcal{F}\big(\Lambda^{-2}\na \D g_n\big)\big|,
\end{align*}
which in turn gives
\bbal
\|t\pa_tu_{2,n}^{\rm{ap}}\|_{H^{s+k}}
&=\Big[\int_{\R^2}(1+|\xi|^2)^{\frac{s+k}{2}}\big|\mathcal{F}\big(t\pa_tu_{2,n}^{\rm{ap}}\big)\big|^2\dd\xi\Big]^{\fr12}\\
&\les\Big[\int_{\mathcal{C}_n }(1+|\xi|^2)^{\frac{s+k}{2}}\big|\mathcal{F}\big(\Lambda^{-2}\nabla f_n\big)\big|^2\dd\xi\Big]^{\fr12}
+t\Big[\int_{\mathcal{B}(0,1)}\big|\mathcal{F}\big(\Lambda^{-2}\na \D g_n\big)\big|^2\dd\xi\Big]^{\fr12}\\
&\les\|f_n\|_{H^{s+k-1}}
+t\|g_n\|_{L^{2}}.
\end{align*}
{\bf A.2\quad Proof of \eqref{a2}.}\quad
From \eqref{a}, then
\bbal
\mathcal{F}\big(\pa_t\varrho_{2,n}^{\rm{ap}}\big)&=
-\frac{\lambda_+\lambda_-(e^{\lambda_+t}-e^{\lambda_-t})}{\lambda_+-\lambda_-}\mathcal{F}\big(f_n\big)
-\frac{\lambda_{+}e^{\lambda_+t}-\lambda_{-}e^{\lambda_-t}}{\lambda_+-\lambda_-}\mathcal{F}\big( \D g_n\big)\\&=
-\frac{e^{\lambda_+t}-e^{\lambda_-t}}{\sqrt{\mu^2-|\xi|^{-2}}}\mathcal{F}\big(f_n\big)
-\frac{\lambda_{+}e^{\lambda_+t}-\lambda_{-}e^{\lambda_-t}}{\lambda_+-\lambda_-}\mathcal{F}\big( \D g_n\big),
\end{align*}
similarly, we have
\bbal
\|\pa_t\varrho_{2,n}^{\rm{ap}}\|_{H^{s+k}}
&=\Big[\int_{\R^2}(1+|\xi|^2)^{\frac{s+k}{2}}\big|\mathcal{F}\big(\pa_t\varrho_{2,n}^{\rm{ap}}\big)\big|^2\dd\xi\Big]^{\fr12}\\
&\les\Big[\int_{\mathcal{C}_n }(1+|\xi|^2)^{\frac{s+k}{2}}\big|\widehat{f_n}(\xi)\big|^2\dd\xi\Big]^{\fr12}
+\Big[\int_{\mathcal{B}(0,1)}(1+|\xi|^2)^{\frac{s+k}{2}}\big|\mathcal{F}\big( \D g_n\big)\big|^2\dd\xi\Big]^{\fr12}\\
&\les\|f_n\|_{H^{s+k}}
+\|g_n\|_{L^{2}}.
\end{align*}
{\bf A.3\quad Proof of \eqref{E}.}\quad Notice that
\bbal
&\widehat{\varrho^{\rm{ap}}_{2,n}}-\widehat{f_n}=\frac{\lambda_+(e^{\lambda_-t}-e^{\lambda_+t})}{\lambda_+-\lambda_-}\widehat{f_n}+(e^{\lambda_+t}-1)\widehat{f_n}
-\frac{e^{\lambda_+t}-e^{\lambda_-t}}{\lambda_+-\lambda_-}\mathcal{F}(\D g_n),\\
&\widehat{u^{\rm{ap}}_{2,n}}-\widehat{g_n}=-\frac{e^{\lambda_+t}-e^{\lambda_-t}}{\lambda_+-\lambda_-}\widehat{\nabla f_n}
+\Big(e^{\lambda_+t}-1+\frac{\lambda_-(e^{\lambda_+t}-e^{\lambda_-t})}{\lambda_+-\lambda_-}\Big)\widehat{g_n},
\end{align*}
then by \eqref{lambda-es}
\bbal
\|\varrho^{\rm{ap}}_{2,n}-f_n\|_{H^{s+1}}
&\les t\Big[\int_{\mathcal{C}_n }(1+|\xi|^2)^{\frac{s+1}{2}}|\widehat{f_n}|^2\dd\xi
+\int_{\mathcal{B}(0,1)}(1+|\xi|^2)^{\frac{s+1}{2}}\big|\mathcal{F}\big(\D g_n\big)\big|^2\dd\xi\Big]^{\fr12}\\
&\les t\big(\|f_n\|_{H^{s+1}}
+\|g_n\|_{L^{2}}\big)\les t2^{n}
\end{align*}
and
\bbal
\|u^{\rm{ap}}_{2,n}-g_n\|_{H^{s+k}}
&\les \Big[\int_{\mathcal{C}_n }(1+|\xi|^2)^{\frac{s+k}{2}}\big|\widehat{\Lambda^{-1}f_n}\big|^2\dd\xi
+t^2\int_{\mathcal{B}(0,1)}(1+|\xi|^2)^{\frac{s+k}{2}}\big|\mathcal{F}\big( g_n\big)\big|^2\dd\xi\Big]^{\fr12}\\
&\les \|f_n\|_{H^{s+k-1}}
+t\|g_n\|_{L^{2}}\les 2^{(k-1)n}+t2^{-n}.
\end{align*}
By Lemma \ref{le3}, we deduce
\bbal
&\|(u_{2,n}^{\rm{ap}}-g_n)_1\pa_1\varrho_{2,n}^{\rm{ap}}\|_{H^{s}}\les\|u_{2,n}^{\rm{ap}}-g_n\|_{H^{s-1}}\|\varrho_{2,n}^{\rm{ap}}\|_{H^{s+1}}+\|u_{2,n}^{\rm{ap}}-g_n\|_{H^{s}}\|\varrho_{2,n}^{\rm{ap}}\|_{H^{s}}\les t+2^{-n},\\
&\|(g_n)_1\pa_1(\varrho_{2,n}^{\rm{ap}}-f_n)\|_{H^{s}}\les\|g_n\|_{H^{s}}\|\varrho_{2,n}^{\rm{ap}}-f_n\|_{H^{s+1}}\les t,\\
&\|(u_{2,n}^{\rm{ap}})_2\pa_2\varrho_{2,n}^{\rm{ap}}\|_{H^{s}}\les\|u_{2,n}^{\rm{ap}}\|_{H^{s-1}}\|\pa_2\varrho_{2,n}^{\rm{ap}}\|_{H^{s}}
+\|u_{2,n}^{\rm{ap}}\|_{H^{s}}\|\pa_2\varrho_{2,n}^{\rm{ap}}\|_{H^{s-1}}\les 2^{-n},
\end{align*}
which implies \eqref{E}.\\
{\bf A.4\quad Proof of \eqref{zw}.}\quad Notice that the support condition of $\widehat{\phi}$, we deduce that
\bal\label{rm1}
2^{-n}\|\pa_1\big(\phi(x_1)\phi(x_2)\big)\pa_1f_n\|_{H^s}\gtrsim&~\|\phi'(x_1)\phi(x_1)\phi^2(x_2)\cos(2^nx_1)\|_{L^2}\nonumber\\&-2^{-n}
\|(\phi'(x_1)\phi(x_2))^2\sin(2^nx_1)\|_{L^2}\nonumber\\
\gtrsim&~\|\phi^2(x_2)\|_{L^2}\|\phi'(x_1)\phi(x_1)\cos(2^nx_1)\|_{L^2}\nonumber\\&-2^{-n}\|(\phi'(x_1)\phi(x_2))^2\|_{L^2}.
\end{align}
Using the simple formula
$$2\cos^2(2^nx_1)=1-\cos(2^{n+1}x_1),$$
we have
$$\|\phi'(x_1)\phi(x_1)\cos(2^nx_1)\|^2_{L^2}=\fr12\|\phi'(x_1)\phi(x_1)\|^2_{L^2}-\fr12\int_{\R}|\phi'(x_1)\phi(x_1)|^2\cos(2^{n+1}x_1)\dd x_1,$$
which follows from Riemann-Lebesgue's Theorem that
\bal\label{rm2}\lim_{n\to\infty}\|\phi'(x_1)\phi(x_1)\cos(2^nx_1)\|^2_{L^2}=\fr12\|\phi'(x_1)\phi(x_1)\|^2_{L^2}.\end{align}
Combining \eqref{rm1} and \eqref{rm2} yields the desired \eqref{zw}.
\section*{Acknowledgments}
J. Li is supported by the National Natural Science Foundation of China (Grant No.11801090). Y. Yu is supported by the Natural Science Foundation of Anhui Province (No.1908085QA05). W. Zhu is partially supported by the National Natural Science Foundation of China (Grant No.11901092) and Natural Science Foundation of Guangdong Province (No.2017A030310634).

\end{document}